\documentclass[a4paper]{amsart}
\usepackage{amsmath,amstext,amsthm,amscd,amsopn,verbatim,amssymb, amsfonts}
\usepackage[bbgreekl]{mathbbol}

\usepackage{fullpage}

\usepackage{tikz}
\usetikzlibrary{matrix}
\usetikzlibrary{shapes}
\usetikzlibrary{arrows}
\usetikzlibrary{calc,3d}
\usetikzlibrary{decorations,decorations.pathmorphing}
\usetikzlibrary{through}
\tikzset{ext/.style={circle, draw,inner sep=1pt},int/.style={circle,draw,fill,inner sep=0, minimum size=5},nil/.style={inner sep=1pt}}
\tikzset{xst/.style={draw, cross out, minimum size=5, }}
\tikzset{exte/.style={circle, draw,inner sep=3pt},inte/.style={circle,draw,fill,inner sep=3pt}}
\tikzset{diagram/.style={matrix of math nodes, row sep=3em, column sep=2.5em, text height=1.5ex, text depth=0.25ex}}
\tikzset{diagram2/.style={matrix of math nodes, row sep=0.5em, column sep=0.5em, text height=1.5ex, text depth=0.25ex}}

\theoremstyle{plain}
  \newtheorem{thm}{Theorem}
  
  \newtheorem{prop}{Proposition}

\theoremstyle{definition}
  \newtheorem{ex}{Example}
  \newtheorem{rem}{Remark}

\newcommand{\co}[2]{\left[{#1},{#2}\right]} % commutator

\newcommand{\R}{{\mathbb{R}}}

\newcommand{\Graphs}{{\mathsf{Graphs}}}

\newcommand{\La}{\Lambda} 

\newcommand{\dGraphs}{{\mathsf{dGraphs}}}
\newcommand{\fdGraphs}{{\mathsf{dfGraphs}}}

\newcommand{\BVGraphs}{{\mathsf{BVGraphs}}}
\newcommand{\fBVGraphs}{{\mathsf{fBVGraphs}}}
\newcommand{\BV}{{\mathsf{BV}}}

\newcommand{\Gra}{{\mathsf{Gra}}}
\newcommand{\EGra}{{\mathsf{EGra}}}

\newcommand{\dGra}{{\mathsf{dGra}}}

\newcommand{\SGra}{{\mathsf{SGra}}}

\newcommand{\bbo}{ \mathbb{1} }

\newcommand{\whl}{\circlearrowleft }

\newcommand{\op}{\mathcal}

\newcommand{\Lie}{\mathsf{Lie}}
\newcommand{\ELie}{\mathsf{ELie}}

\newcommand{\Def}{\mathrm{Def}}

\newcommand{\GC}{\mathsf{GC}}
\newcommand{\EGC}{\mathsf{EGC}}
\newcommand{\fEGC}{\mathsf{fEGC}}
\newcommand{\fGC}{\mathsf{fGC}}
\newcommand{\fcGC}{\mathsf{fcGC}}

\newcommand{\dfGC}{\mathsf{dfGC}}
\newcommand{\fdGC}{\mathsf{dfGC}}

\newcommand{\vout}{\mathit{out}}
\newcommand{\vin}{\mathit{in}}

\newcommand{\bpm}{\begin{pmatrix}}
\newcommand{\epm}{\end{pmatrix}}
\newcommand{\Tpoly}{T_{\rm poly}}

\newcommand{\Dpoly}{D_{\rm poly}}

\newcommand{\grt}{{\mathfrak{grt}_1}}

\newcommand{\mU}{\mathcal{U}}
\newcommand{\mV}{\mathcal{V}}

\DeclareMathOperator{\gr}{gr}
\DeclareMathOperator{\End}{End}

\newcommand{\gra}{\mathrm{gra}}
\newcommand{\dgra}{\mathrm{dgra}}
\newcommand{\K}{\mathbb{K}}
\DeclareMathOperator{\Tw}{\mathrm{Tw}}

\begin{document}
\title{The Grothendieck-Teichm\"uller group action on differential forms and formality morphism of chains}
\author{Thomas Willwacher}
\address{Institute of Mathematics\\ University of Zurich\\ Winterthurerstrasse 190 \\ 8057 Zurich, Switzerland}
\email{thomas.willwacher@math.uzh.ch}

%\thanks{The author was partially supported by the Swiss National Science Foundation (grant 200020-105450).}
% \subjclass[2000]{16E45; 53D55; 53C15; 18G55}
% \date{}
\keywords{Formality, Deformation Quantization, Graph Cohomology}

\begin{abstract}
It is known that one can associate a Kontsevich-type formality morphism to every Drinfeld associator. We show that this morphism may be extended to a Kontsevich-Shoikhet formality morphism of cochains and chains, by describing the action of the Grothendieck-Teichm\"uller group on such objects (up to homotopy).
\end{abstract}
\maketitle

\section{Introduction}
Let $\Tpoly$ be the space of multivector fields on $\R^n$ and let $\Dpoly$ be the space of multidifferential operators on $\R^n$. M. Kontsevich's formality Theorem \cite{K1} states that there is a $\Lie_\infty$ quasi-isomorphism
\[
\mU: \Tpoly[1] \to \Dpoly[1].
\] 
Here we understand $\Tpoly[1]$ as a Lie algebra endowed with the Schouten bracket and $\Dpoly[1]$ as a Lie algebra endowed with the Gerstenhaber bracket.
The differential forms $\Omega_\bullet$ on $\R^n$, with non-positive grading, form a Lie module over $\Tpoly[1]$. The action of a $k$-vector field $\gamma$ on a differential form $\alpha$ is given by the Lie derivative
\[
 L_\gamma \alpha = d \iota_\gamma \alpha + (-1)^{k} \iota_\gamma \alpha
\]
where $d$ is the de Rham differential and $\iota_\gamma$ is the operation of contraction with $\gamma$. Similarly the completed Hochschild chain complex $C_\bullet=C_\bullet(C^\infty(\R^n), C^\infty(\R^n))$ forms a module over the multidifferential operators $\Dpoly$, see \cite{tsygan}. 
B. Shoikhet \cite{shoikhet} showed that there is a $\Lie_\infty$ quasi-isomorphism of modules
\[
\mV \colon C_\bullet \to \Omega_\bullet
\]
thus proving an earlier conjecture of B. Tsygan \cite{tsygan}. Here $C_\bullet$ is considered as a $\Lie_\infty$ module over $\Tpoly[1]$ by pulling back the $\Dpoly[1]$ module structure along the morphism $\mU$.

In fact, the ``correct'' objects to consider are not formality morphisms $\mU$ on $\R^n$ for some fixed $n$, but \emph{stable} formality morphisms in the sense of \cite{vasilystable}. The components of such morphisms are expressed by operations in a suitable operad of graphs, which acts on the pair $(\Tpoly, \Dpoly)$ for any $n$, so that from a stable formality morphism one obtains a formality morphism for each $n$. By the results of V. Dolgushev \cite{vasilystable} the space of stable formality morphisms up to homotopy is a torsor for the zeroth cohomology of the graph complex $\fGC$ (see below for a definition). The latter object may be identified with the Grothendieck-Teichm\"uller Lie algebra $\grt$, see \cite{grt}. It follows that the space of stable formality morphisms (up to homotopy) may be identified with the space of Drinfeld associators. An explicit construction of a formality morphism given a Drinfeld associator has been described by D. Tamarkin earlier \cite{tamanother,hinich}. It was shown in \cite{brformality} that applying one version of D. Tamarkin's construction to the Alekseev-Torossian associator \cite{ATassoc,pavol} yields Kontsevich's formality morphism. Furthermore all formality morphisms thus obtained may actually be extended to homotopy Gerstenhaber (instead of just $\Lie_\infty$) formality morphisms \cite{brformality}.

In this paper we extend the above picture to formality morphisms of chains and cochains. Concretely, we will show that one can associate a formality morphism of chains and cochains to every Drinfeld associator.

\begin{thm}\label{thm:cor}
For each Drinfeld associator $\Phi$ there is a stable formality morphism $\mU_\Phi$ in the homotopy class associated to $\Phi$, together with a formality morphism of chains 
\[
 \mV_\Phi : C_\bullet \to \Omega_\bullet
\]
also given by graphical formulas, where the $\Tpoly[1]$ action on $C_\bullet$ is obtained by pulling back the action of $\Dpoly[1]$ on $C_\bullet$ along $\mU_\Phi$.
\end{thm}

To show the Theorem it suffices to lift the action of the Grothendieck-Teichm\"uller group on the homotopy classes of stable formality morphisms of cochains to homotopy classes of stable formality morphisms of cochains and chains. 
% 
% Hence the most interesting questions in this regard have been settled, and one has a pretty complete picture.
% However, one also wants to extend this picture to include formality morphisms on chains as well.
% A part of the considerations above has an ``obvious'' analogue in the chains setting. For example, there is a natural notion of stable formality morphism of chains and cochains. 
To this end we will consider a version of the graph complex which we call $\fEGC$ acting on $(\Tpoly[1], \Omega_\bullet)$ by $\Lie_\infty$ derivations, and hence on (stable) formality morphisms of chains and cochains. 
%One may use the methods of \cite{vasilystable} or \cite{mestable} to show that the induced action of $H^0(\fEGC)$ is free and transitive on homotopy classes of stable formality morphisms of chains and cochains, though this will not be done here.
%I claim that one may also show that each such formality morphism of chains and cochains may be extended to a homotopy Kontsevich-Soibelman morphism, using the fact that $\fEGC$ acts on a graphical model of the Kontsevich-Soibelman operad, see \cite{brformality}. This will also not be done here.
%Concerning this graph complex, there is however a small gap in our understanding, whose solution needs a slightly new idea beyond an ``extension to the chains case'' of ideas already existing in the literature.
More concretely, as a graded lie algebra
\[
 \fEGC \cong \fGC \ltimes \fGC_1
\]
where $\fGC$ is as before and acts on $\Tpoly[1]$ while the part $\fGC_1$ to be introduced below acts on $\Omega_\bullet$ considered as $\Lie_\infty$ module.
There is a 
%Namely, one has to compute the cohomology, in particular in degree zero, of the Lie algebra $\fEGC$ governing stable formality morphism of chains and cochains. %That is the purpose of this paper.
%However, it is not clear a priori what the relation of $\GC$ and $\EGC$ is. 
projection map $\fEGC\to \fGC$, but a priori it is not clear that, for example, any cocycle in $\fGC$ may be extended to one in $\fEGC$.
To describe the cohomology of $\fEGC$, we need to introduce the ``divergence'' operator $\nabla$ on $\fGC$, which is defined on a graph by summing over all ways to add an edge. Alternatively,
\[
 \nabla = [
 \begin{tikzpicture}[baseline=-.65ex]
\clip (-.3,-.5) rectangle (.3, .5);  
  \node[int] (v) at (0,0) {};
\draw (v) edge[loop, looseness=10] (v);
 \end{tikzpicture}
,\cdot ]
\]
is the Lie bracket with a tadpole (short-loop) graph. The operator $\nabla$ commutes with the differential and induces a degree $-1$ operator on $H(\fGC)$. Unfortunately, the precise form of this operator on cohomology is unknown, owed to the fact that most of the cohomology of $\fGC$ is unknown.
The main result of this paper is the computation of the cohomology of $\fEGC$ in terms of that of $\fGC$, and the action of $\nabla$.
\begin{thm}
\label{thm:main}
\[
 H(\fEGC) \cong \K B \oplus \K B \otimes H(\fGC) \oplus \K \bbo \oplus H(\fGC\oplus \fGC[-2], \nabla)
\]
where $B$ and $\bbo$ are explicitly known cohomology classes described below (see \eqref{equ:oneandB}), and the $\nabla$ on the right is understood as a degree 1 map from $\fGC$ to $\fGC[-2]$.
 Furthermore, the map $H^0(\EGC)\to H^0(\GC)\cong \grt$ is an isomorphism.
\end{thm}
%In fact, we show a bit more. We will give a complete description of the cohomlogy of $\EGC$ below. 
We provide an explicit combinatorial formula for the cocycles in $\fEGC$ corresponding to ``divergence free'' (see section \ref{sec:divergence} below) cocycles in $\fGC$.

Hence the Grothendieck-Teichm\"uller group action on the homotopy classes of stable formality morphisms of cochains lifts to an action on (stable) formality morphisms of cochains and chains. In particular, we may associate a stable formality morphism of chains and cochains to each Drinfeld associator $\Phi$. Concretely, to the Alekseev-Torossian associator we associate the Kontsevich-Shoikhet morphism. The stable formality morphisms corresponding to another Drinfeld associators $\Phi'$ may be recovered by acting with the unique element of the Grothendieck-Teichm\"uller group sending $\Phi$ to $\Phi'$.
This shows Theorem \ref{thm:cor}.

%Note that by the results of \cite{brformality} the  morphism to a homotopy Kontsevich-Soibelman, or a $\calc_\infty$ morphism, cf. also \cite{DTT}. 
%The stable formality morphisms corresponding to other Drinfeld associators may be recovered by the $\grt\cong H^0(\EGC)$ action. 

 %\item All such morphisms may in fact be extended to $\calc_\infty$ morphisms, see \cite{brformality}.
%\end{itemize}

 %(TODO: WHAT WE SHOULD ALSO SHOW IS THAT THIS EXHAUSTS ALL STABLE CHAINS/COCHAINS FORMALITY MORPHISMS. I AM SURE THIS IS TRUE SINCE EGC GOVERNS STABLE AUTOMORPHISMS OF $(\Tpoly,\Omega)$, AND STABLE FORMALITY MORPHISMS ARE JUST STABLE QISOMORPHISMS TO $(\Dpoly, C)$. BUT ONE SHOULD VERIFY THIS... WOULD BE A GOOD PROJECT FOR A STUDENT.)

\subsection{Structure of the paper}
In section \ref{sec:graphcomplexes} we briefly recall the definition of the graph complex $\fGC$, and introduce the graph complex $\fEGC$, which can be extracted from \cite{brformality}.
Section \ref{sec:theproof} contains the proof of Theorem \ref{thm:main}.
Finally, in section \ref{sec:divexplicit} we will derive the explicit formula for the cocycles in $\fEGC$ corresponding to divergence free cocyles in $\fGC$.

\subsection*{Acknowledgements}
The author was partially supported by the Swiss National Science foundation, grants PDAMP2\_137151 and 200021\_150012.

\section{Graph complexes and operads}\label{sec:graphcomplexes}

\subsection{Definition of the graph complexes}
Let $\gra_{N,k}$ be the set of undirected graphs with vertex set $[N]=\{1,\dots, N\}$ and edge set $[k]$. It carries an action of the group $S_N\times S_k$ by renumbering the vertices and renumbering the edges.

Fix a field $\K$ of characteristic zero. One may define an operad of graphs $\Gra$ such that the space of $N$-ary operations is
\[
 \Gra(N) := \prod_{k\geq 0} \left( \K\langle \gra_{N,k} \rangle \otimes \K[1]^{\otimes k} \right)_{S_k},
\]
where $S_k$ acts diagonally on $\gra_{N,k}$ and on the factors of $\K[1]$ by permutations, with appropriate signs.
Elements of $\Gra$ are series of undirected graphs with edges of degree $-1$, for example the following one:
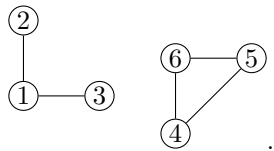
\[
 \begin{tikzpicture}
  \node[ext] (v1) at (0,0) {1};
\node[ext] (v2) at (0,1) {2};
\node[ext] (v3) at (1,0) {3};
\node[ext] (v4) at (2,-0.5) {4};
\node[ext] (v5) at (3,.5) {5};
\node[ext] (v6) at (2,.5) {6};
\draw (v1) edge (v2) edge (v3)
(v4) edge (v5) edge (v6) (v5) edge (v6);
 \end{tikzpicture}.
\]
The operad structure is given by inserting one graph into a vertex of another and reconnecting the incoming edges in all possible ways, see the introductory sections of \cite{grt}.

Analogously, let $\dgra_{N,k}$ be the set of directed graphs with vertex set $[N]=\{1,\dots, N\}$ and edge set $[k]$, and define an operad $\dGra$ whose space of $N$-ary operations is
\[
 \dGra(N) := \prod_{k\geq 0} \left( \K\langle \dgra_{N,k} \rangle \otimes \K[1]^{\otimes k} \right)_{S_k}.
\]
There is a map of operads 
\[
\Gra\to \dGra 
\]
sending each undirected edge to the sum of the edges in either direction.

Next consider a similar set $\dgra_{1,N,k}$ whose elements are directed graphs with vertex set $N\cup\{\vin,\vout\}$ and edge set $[k]$, such that at the vertex $\vout$ there are only outgoing edges and at the vertex $\vin$ there are only incoming edges.
We define vector spaces 
\[
 \Gra_1(N) := \prod_{k\geq 0} \left( \K\langle \gra_{1,N,k} \rangle \otimes \K[1]^{\otimes k} \right)_{S_k}.
\]
Elements are linear combinations of directed graphs with two special vertices $\vin$ and $\vout$, for example the following:
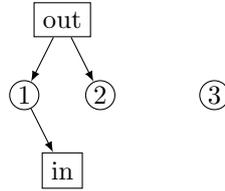
\[
 \begin{tikzpicture}[every edge/.style={draw, -latex}, baseline=-.65ex]
 \node[draw] (out) at (0,1) {out}; 
 \node[draw] (in) at (0,-1) {in};
 \node[ext] (v1) at (-.5,0) {1};
 \node[ext] (v2) at (.5,0) {2};
 \node[ext] (v3) at (2,0) {3};
\draw (out) edge (v1) edge (v2) (v1) edge (in);
\end{tikzpicture}
\]
The spaces $\Gra(N)$ and $\Gra_1(N)$ assemble into a 2-colored operad $\EGra$. Here $\Gra(N)$ is the space of operations with $N$ inputs and the output in color 1, and $\Gra_1(N)$ is the space of operations with $N$ inputs in color 1, one input in color $2$, represented by the vertex $\vin$, and the output in color 2, represented by the vertex $\vout$.
The operadic compositions are obtained by inserting graphs at vertices of others and reconnecting the dangling edges in all possible ways, see \cite{brformality}.
Denote the Lie operad by $\Lie$, its suspension by $\La\Lie$, and the minimal cofibrant resolution thereof by $\La\Lie_\infty$.
We denote the two-colored operad that governs a $\La\Lie$ algebra and a module by $\ELie$, and its minimal cofibrant resolution by $\ELie_\infty$.

Recall from loc. cit. that there is a natural map of colored operads $\ELie_\infty\to \ELie\to \EGra$.
We define the full graph complex as the operadic deformation complex
\[
 \fEGC := \Def(\ELie_\infty\to \EGra).
\]
Concretely, elements of $\fEGC$ are just series of graphs as occuring in $\Gra$ and $\SGra$, invariant under permutations of the symbols $1,2,\dots$ decorating the vertices.
The abstract definition as a deformation complex has the advantage that it is immediate that one has a differential graded (dg) Lie algebra structure on $\fEGC$.

Let us disect $\fEGC$ again into smaller pieces. First it contains a quotient Lie algebra
\[
 \fGC := \Def(\La\Lie_\infty\to \Gra),
\]
the full graph complex. It consists of series of graphs as in $\Gra$, invariant under permutations of the symbols $1,2,\dots$ decorating the vertices.

Similarly, there is a sub-dg Lie algebra $\fGC_1 \subset \fEGC$ consisting of series in graphs as in $\Gra_1$, invariant under renumbering the vertices.
One can write $\fEGC = \fGC \ltimes \fGC_1$.

A similar construction can be carried out using the directed graphs operad $\dGra$. One in particular obtains the directed graph complex 
\[
\dfGC \cong \Def(\La\Lie_\infty\to \Gra).
\]
It can be checked (see \cite[Appendix K]{grt}) that the natural map $\fGC\to \dfGC$ is a quasi-isomorphism. 

\begin{rem}
 For concreteness, let us give pictorial description of the differential on $\fEGC$.
It is a sum of three pieces, the first piece $\delta:\fGC\to \fGC$ sums over all vertices and splits the vertex into two in all possible ways. Pictorially
\[
 \delta 
\begin{tikzpicture}[baseline=-.65ex]
 \node[int] (v) at (0,0){};
\draw (v) edge +(-.5,.5) edge +(0,.5) edge +(.5,.5) edge +(-.5,-.5) edge +(0,-.5) edge +(.5,-.5);
\end{tikzpicture}
=
 \sum 
\begin{tikzpicture}[baseline=-.65ex]
 \node[int] (v) at (0,-0.3){};
\node[int] (w) at (0,0.3){};
\draw (w) edge +(-.5,.5) edge +(0,.5) edge +(.5,.5) (v) edge +(-.5,-.5) edge +(0,-.5) edge +(.5,-.5) edge (w);
\end{tikzpicture}.
\]
The second piece which we also denote by $\delta$ maps $\fGC_1$ to $\fGC_1$ as follows:
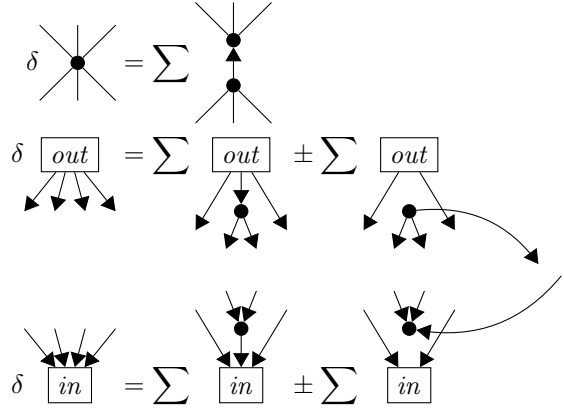
\begin{align*}
  \delta 
\begin{tikzpicture}[baseline=-.65ex]
 \node[int] (v) at (0,0){};
\draw (v) edge +(-.5,.5) edge +(0,.5) edge +(.5,.5) edge +(-.5,-.5) edge +(0,-.5) edge +(.5,-.5);
\end{tikzpicture}
&=
 \sum 
\begin{tikzpicture}[baseline=-.65ex]
 \node[int] (v) at (0,-0.3){};
\node[int] (w) at (0,0.3){};
\draw (w) edge +(-.5,.5) edge +(0,.5) edge +(.5,.5) (v) edge +(-.5,-.5) edge +(0,-.5) edge +(.5,-.5) edge[-triangle 60] (w);
\end{tikzpicture}
\\ 
\delta 
 \begin{tikzpicture}[baseline=-.65ex]
\node[draw] (v) at (0,0) {$\vout$};
\draw[-triangle 60] (v) edge +(-.6,-.75) edge +(-.2,-.75) edge +(.2,-.75) edge +(.6,-.75);
\end{tikzpicture}
&=\sum
\begin{tikzpicture}[baseline=-.65ex]
\node[draw] (v) at (0,0) {$\vout$};
 \node[int] (w) at (0,-.75){};
 \draw[-triangle 60] (w) edge +(-.2,-.5) edge +(0.2,-.5) (v) edge +(-.6,-1) edge +(.6,-1) edge (w);
\end{tikzpicture}
\pm
\sum
\begin{tikzpicture}[baseline=-.65ex]
\node[draw] (v) at (0,0) {$\vout$};
 \node[int] (w) at (0,-.75){};
 \draw[-triangle 60] (w) edge +(-.2,-.5) edge +(0.2,-.5) (v) edge +(-.6,-1) edge +(.6,-1);
 \draw[-triangle 60] (w) edge[bend left] +(1.7,-.7);
\end{tikzpicture}
\\
\delta 
 \begin{tikzpicture}[baseline=-.65ex]
\node[draw] (v) at (0,0) {$\vin$};
\draw[triangle 60-] (v) edge +(-.6,.75) edge +(-.2,.75) edge +(.2,.75) edge +(.6,.75);
\end{tikzpicture}
&=\sum
\begin{tikzpicture}[baseline=-.65ex]
\node[draw] (v) at (0,0) {$\vin$};
 \node[int] (w) at (0,.75){};
 \draw[triangle 60-] (w) edge +(-.2,.5) edge +(0.2,.5) (v) edge +(-.6,1) edge +(.6,1) edge (w);
\end{tikzpicture}
\pm
\sum
\begin{tikzpicture}[baseline=-.65ex]
\node[draw] (v) at (0,0) {$\vin$};
 \node[int] (w) at (0,.75){};
 \draw[triangle 60-] (w) edge +(-.2,.5) edge +(0.2,.5) (v) edge +(-.6,1) edge +(.6,1);
 \draw[triangle 60-] (w) edge[bend right] +(2,.7);
\end{tikzpicture}
\end{align*}
Here the bent edges are supposed to be reconnected to some other vertex of the graph. (One sums over all choices.)

Finally, there is a third piece of the differential $\delta_1\colon \fGC\to\fGC_1$, which sends an element $\Gamma\in \fGC$ to 
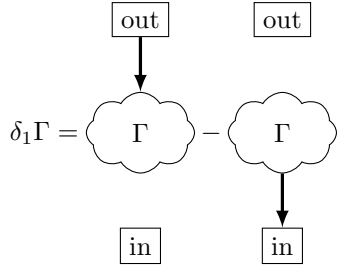
\[
\delta_1\Gamma =
 \begin{tikzpicture}[every edge/.style={draw, very thick, -latex}, baseline=-.65ex]
\node [cloud, draw,cloud puffs=8,cloud puff arc=120, aspect=1.6, inner ysep=.5em] (v1) at (0,0) {$\Gamma$};
\node[draw] (out) at (0,1.5) {out}; 
 \node[draw] (in) at (0,-1.5) {in};
\draw (out) edge (v1); 
\end{tikzpicture}
-
\begin{tikzpicture}[every edge/.style={draw, very thick, -latex}, baseline=-.65ex]
\node [cloud, draw,cloud puffs=8,cloud puff arc=120, aspect=1.6, inner ysep=.5em] (v1) at (0,0) {$\Gamma$};
\node[draw] (out) at (0,1.5) {out}; 
 \node[draw] (in) at (0,-1.5) {in};
\draw (v1) edge (in);.
\end{tikzpicture}
\]
\end{rem}

\begin{rem}
The graph complexes above admit disconnected graphs and vertices of all valences. One often restricts to smaller subcomplexes. For example, there are further dg Lie subalgebras
\[
 \GC \subset \fcGC \subset \fGC
\]
where $\fcGC$ consists of the connected graphs only and $\GC$ consists of connected graphs with all vertices of valence $\geq 3$.
Clearly,  $\fGC\cong S^+(\fcGC[-2])[2]$ is a completed symmetric product without unit.
Furthermore, the cohomology of $\fcGC$ may be expressed in terms of that of $\GC$, cf. \cite[Proposition 3.4]{grt}.
Similarly, we may identify sub-dg Lie algebras
\[
 \GC_1 \subset \fcGC_1 \subset \fGC_1.
\]
Here $\fcGC_1$ consists of graphs that are connected after deleting vertices $\vin$ and $\vout$. Such graphs will be called \emph{internally connected}. The dg Lie algebra  
$\GC_1$ consists of graphs that are internally connected and all internal vertices, i.~e., vertices other than $\vin$ or $\vout$ are at least two-valent.
We finally define the dg Lie subalgebra 
\[
 \EGC := \GC \ltimes \GC_1 \subset  \fEGC.
\]
\end{rem}

%There are two important subcomplexes $\GC_n \subset \fGCc_n\subset \fGC_n$. Here $\fGCc_n$ consists of series of connected graphs, and $\GC_n$ consists of series in connected graphs all of whose vertices are at least trivalent.
%The cohomology of $\fGC_n$ is clearly just the symmetric product of that of $\fGCc_n$. Furthermore, the latter cohomology is the sum of of $\GC_n$ and the wheel classes, which are represented by graphs which are loops of bivalent vertices.

% \begin{rem}
%  In this paper we take the approach that all graphs are allowed to contain short cycles a priori, i.~e., edges connecting some vertex to itself.
% This is not consistent with the notation used elsewhere. However, we don't want to spoil the already complicated notation further by putting ${}^\circlearrowleft$ superscripts everywhere. 
% \end{rem}

% It is a sub-dg Lie algebra of the deformation complex
% \[
% \Def(\La\Lie_\infty\to \Gra). 
% \]
Since $\Gra$ (and $\dGra$) can be naturally represented on $\Tpoly$, degree zero cocycles in $\fGC$ act on $\Tpoly[1]$ by $\Lie_\infty$ derivations, and hence also on the space of formality morphisms $\Tpoly[1]\to \Dpoly[1]$. Similarly, the operad $\EGra$ may be naturally represented on the pair (i.~e., on the colored vector space) $(\Tpoly,\Omega_\bullet)$, and hence degree zero cocycles in $\fEGC$ act naturally on the pair $(\Tpoly[1],\Omega_\bullet)$ by $\Lie_\infty$ derivations and $\Lie_\infty$ derivations of modules.
Concretely, let $x+x_1$ be a closed degree zero element of $\fEGC$, with $x\in \fGC$, $x_1\in \fGC_1$. Then $x$ acts, as before, by a $\Lie_\infty$ derivation on $\Tpoly[1]$. However, $x$ does not necessarily respect the $\Lie_\infty$ module structure on $\Omega_\bullet$, but changes the $\Lie_\infty$ module structure to an ``infinitesimally different'' one. The element $x_1$ corrects for this defect by providing an ``infinitesimal $\Lie_\infty$ morphism of modules'' between $\Omega_\bullet$ and $\Omega_\bullet$ with the modified module structure. 

It was shown in \cite{grt} that $H^0(\GC)\cong \grt$, and hence one recovers the action of the Grothendieck-Teichm\"uller group on formality morphisms of cochains.
The main technical contribution of this work is to extend this action to the formality morphisms on chains $C_\bullet\to \Omega_\bullet$ by showing that $H^0(\fEGC)\cong \grt$.

\subsection{Some cocycles in $\fGC_1$}
Let us decribe some cocyles in $\fGC_1$. The simplest two are the following
\begin{align}\label{equ:oneandB}
 \bbo &=
\begin{tikzpicture}[every edge/.style={draw, -latex}, baseline=-.65ex]
 \node[draw] (out) at (0,.7) {out}; 
 \node[draw] (in) at (0,-.7) {in};
\end{tikzpicture}
&
 B=
\begin{tikzpicture}[every edge/.style={draw, -latex}, baseline=-.65ex]
 \node[draw] (out) at (0,.7) {out}; 
 \node[draw] (in) at (0,-.7) {in};
 \draw (out) edge (in);.
\end{tikzpicture}.
\end{align}

Furthermore, let us note that there is a map of complexes $\dfGC\to \fGC_1$ mapping a graph $\Gamma$ to $I\bullet \Gamma$, where $I$ is the graph with one internal vertex and no edges, i. e.,
\[
 I=\,
\begin{tikzpicture}[every edge/.style={draw, -latex}, baseline=-.65ex]
 \node[draw] (out) at (0,.7) {out}; 
 \node[draw] (in) at (0,-.7) {in};
 \node[int] (v1) at (0,0) {};
\end{tikzpicture}
\]
and the $\bullet$ denotes insertion of $\Gamma$ in place of the black vertex.
Note that we also have a map $\fGC\to\fGC_1$ by composition of the previous map with the embedding $\fGC\to\dfGC$.
%, and that this map restricts to a map $\GC\to \GC_1$.

% Note also that there are closed elements $W_k$ in $\fGC$ given by the loop graphs
% \[
%  W_k=
%  \begin{tikzpicture}[baseline=-.65ex]
%   \node (v0) at (0:1) {$\cdots$};
% \node[int] (v1) at (60:1) {};
% \node[int] (v2) at (120:1) {};
% \node[int] (v3) at (180:1) {};
% \node[int] (v4) at (240:1) {};
% \node[int] (v5) at (300:1) {};
% \draw (v0) edge (v1) edge (v5) (v2) edge (v1) edge (v3) (v4) edge (v3) edge (v5);
%  \end{tikzpicture}
% \text{   ($k=5,9,11,\dots$ vertices)}.
% \]
% Recall the following result from \cite{grt}.
% \begin{prop}[see Proposition 3.4 of \cite{grt}]\label{prop:GCfGC}
%  \[
%   H(\fcGC)\cong H(\GC) \oplus \prod_{k=5,9,11,\dots} \K W_k.
%  \]
% \end{prop}

% The cocyles $W_k$ give rise to cocyles in $\GC_1$ by the map above. %elements under the map $\Gamma\mapsto I\bullet \Gamma$ encountered above.

% \begin{rem}
%  The colored operad $\EGra$ acts naturally on the colored vector space $\Tpoly\oplus \Omega_\bullet$. It follows that one obtains an action of $\EGC$. Concretely, let $x+x_1$ be a closed degree zero element of $\EGC$, with $x\in \GC$, $x_1\in \GC_1$. Then $x$ acts, as before, by a $\Lie_\infty$ derivation on $\Tpoly$. However, $x$ does not necessarily respect the $\Lie_\infty$ module structure on $\Omega_\bullet$, but changes the $\Lie_\infty$ module structure to an ``infinitesimally different'' one. The element $x_1$ corrects for this defect by providing an ``infinitesimal $\Lie_\infty$ morphism of modules'' between $\Omega_\bullet$ and $\Omega_\bullet$ with the modified module structure. 
% \end{rem}

\subsection{The divergence operation}\label{sec:divergence}
The graphs occurring in the operads $\Gra$ and $\dGra$ were not allowed to contain tadpoles, i.~e., edges connecting a vertex to itself.
In fact, we might have allowed them as well, to obtain operads we denote by $\Gra^\whl$ and $\dGra^\whl$, and graph complexes $\fGC^\whl$ and $\dfGC^\whl$.
These graph complexes contain their tadpole-free relatives as subcomplexes.
It is shown in Proposition 3.4 of \cite{grt} that the graph complexes with tadpoles are quasi-isomorphic to those without tadpoles except possible for the occurrence of the cohomology class represented by the graph
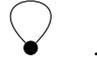
\[
 \begin{tikzpicture}
  \node[int] (v) at (0,0) {};
\draw (v) edge[loop, looseness=20] (v);
 \end{tikzpicture}.
\]
In fact, one may check that the Lie bracket of this graph with itself vanishes, and hence one may define an additional differential of degree -1 on the graph complex $\fGC^\whl$ as 
\[
 \nabla = [
 \begin{tikzpicture}[baseline=-.65ex]
\clip (-.3,-.5) rectangle (.3, .5);  
  \node[int] (v) at (0,0) {};
\draw (v) edge[loop, looseness=10] (v);
 \end{tikzpicture}.
,\cdot].
\]
We call this operator the divergence operator.
It acts on a graph by adding an additional edge in all possible ways.
In fact, the action of $\nabla$ leaves invariant the subspace $\fGC\subset\fGC^\whl$, and hence descends to an operator on that space, that we will also denote by $\nabla$.

In general, it is not known what the induced action of $\nabla$ on the cohomology $H(\fGC)$ is. However, since $H^{<0}(\fGC)=0$ (see \cite{grt}) the subspace $H^0(\fGC)\cong \grt$ is sent to 0 by degree reasons. It follows in particular that the graph cohomology classes corresponding to $\grt$-elements may be represented by divergence-free cocycles, cf. \cite{sergeime}.

\begin{rem}
 Note that in \cite{carlothomas} explicit integral formulas for divergence-free graph cocycles corresponding to all Deligne-Drinfeld elements $\sigma_3, \sigma_5,\dots \in \grt$ are given. % correspond to divergence free graph cocycles, i.~e., elements of $\GC_{closed,div}^0$. Since these elements conjecturally generate $\grt$, it is expected that any element of $\grt$ is represented by a divergence free graph cocycle.
\end{rem}

Below we will need to consider the divergence free sub-dg Lie algebra $\fGC_{div}\subset\fGC$ and similarly $\fdGC_{div}\subset\fdGC$ spanned by the elements of $\fGC$ (respectively of $\fdGC$) closed under $\nabla$, i.~e.,
\[
 \fGC_{div} := \{\gamma\in \fGC \mid \nabla \gamma=0 \}.
\]

% \subsection{Reformulation of the main result}
% Now we may rephrase our main result:
% \begin{thm}
% \label{thm:main}
% The map $\K \bbo \oplus \K B \oplus \fcGC \to \GC_1$ is a quasi-isomorphism.
% Furthermore, there is an explicit map 
% \[
% \GC_{closed,div}^0 \to \EGC
% \]
% where $\GC_{closed,div}^0$ is the subspace of closed, divergence free degree zero elements of $\GC$.
% \end{thm}

% \begin{rem}
% Note that 
% \[
% \GC_{closed,div}^0/ (\delta \GC^{-1} \cap \GC_{closed,div}^0) \cong \grt.
% \]
% Hence the second assertion of Theorem \ref{thm:main} should be seen as providing a somewhat explicit description of the action of $\grt$ in deformation quantization.
% (There is however no known explicit injective ($\Lie_\infty$) map $\grt\to \GC_{closed,div}^0$.)
% \end{rem}

\subsection{The operad $\fdGraphs$ }
We need one more ingredient from the theory of graphical operads, namely a variant of the operad $\Graphs$ introduced by M. Kontsevich in \cite{K2}.
Concretely, consider the operad $\fdGraphs:= \Tw\dGra$ obtained from $\dGra$ by operadic twisting (see \cite{vastwisting}).
Elements of $\Tw\dGra(N)$ are series of directed graphs with $N$ numbered ``external'' vertices $1,\dots, N$ and an arbitrary number of unlabelled ``internal'' vertices.
An example (for $N=1$) can be found in the following picture:
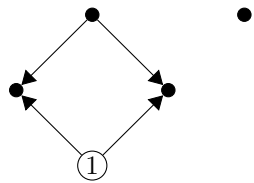
\[
 \begin{tikzpicture}
\node[ext] (v) at (0,0) {1};
\node[int] (v1) at (1,1) {};
\node[int] (v2) at (0,2) {};
\node[int] (v3) at (-1,1) {};
\node[int] (v4) at (2,2) {};
\draw[-triangle 60] (v) edge (v1) edge (v3) (v2) edge (v3) edge (v1);
 \end{tikzpicture}
\]

There is a suboperad $\dGraphs \subset \Tw\dGra$ given by restricting to graphs which do not have connected components with only internal vertices.
For example, the graph above would not have been admissible.
One can check (see \cite{grt}) that $H(\dGraphs)\cong e_2$ is the Gerstenhaber operad.

Similarly, we may consider an operad $\fdGraphs^\whl$, defined in the exactly the same manner, except that the graphs occurring may have tadpoles, i.~e., edges connecting a vertex to itself.
This operad contains an operadic ideal $I\subset \fdGraphs^\whl$ whose elements are series of graphs that contain a tadpole at an internal vertex.
We may form the quotient
\[
 \fBVGraphs := \fdGraphs^\whl / I.
\]
In other words, one sets graphs with tadpoles at internal vertices to zero, while tadpoles at external vertices are still admissible.
The cohomology of the operad $\fBVGraphs$ may be computed, modulo the graph cohomology.

\begin{prop}\label{prop:fBVcohom}
\[
H(\fBVGraphs) \cong \BV \otimes S(H(\fcGC)[-2]).
\]
where $\BV$ is the Batalin-Vilkovisky operad and $S(\dots)$ denotes the completed symmetric product.
\end{prop}
Concretely, the elements of $\BV$ are represented by graphs without internal vertices. The factor $S(H(\fcGC)[-2])$ corresponds to additional connected components fully consisting of internal vertices that may be present.
\begin{proof}[Proof sketch.]
 The graphs contributing to $\fBVGraphs$ are the same as those contributing to $\fdGraphs$, except that there may be one or more tadpoles at some external vertices. 
One checks that these extra tadpoles are not ``seen'' by the differential. Hence the cohomology of $\fBVGraphs$ is the same as that of $\fdGraphs$, with the possible addition of tadpoles at the external vertices.
\end{proof}

\begin{rem}
The operad $\fBVGraphs$ contains a sub-operad $\BVGraphs\subset \fBVGraphs$ formed by graphs which do not have connected components with only internal vertices, that is quasi-isomorphic to the Batalin-Vilkovisky operad.
\end{rem}

An important fact about the operad $\fdGraphs$ is that it is acted upon by the dg Lie algebra $\fdGC$ by operadic derivations.
This follows directly from the formalism of operadic twisting, see \cite{vastwisting}. Let us briefly recall how to obtain the action, leaving the discussion of sign and combinatorial prefactor subleties to loc. cit.
First, there is a right action of $\fdGC$ on $\fdGraphs$ (which is not compatible with the differential) 
\[
 \Gamma \bullet \gamma := \sum_v \Gamma \bullet_v \gamma
\]
where $\gamma\in \fdGC$, $\Gamma\in \fdGraphs$, the sum is over all internal vertices of $\Gamma$, and the notation $\bullet_v$ shall indicate that one inserts $\gamma$ in place of vertex $v$ and reconnects the edges incident to $v$ in all possible ways to vertices of $\gamma$.

Similarly, there is a map of vector spaces $\fdGC\to \fdGraphs(1)$, sending $\gamma \in  \fdGC$ to an element $\gamma_1$ obtained by summing over all vertices of $\gamma$ and declaring the vertex to be external.
In any operad $\op P$ the unary operations $\op P(1)$ form an algebra that acts on the operad by derivations. We denote this action symbolically by
\[
 x \cdot y = x \circ y \pm y\circ x.
\]
Hence we obtain another action of $\fdGC$ on $\fdGraphs$ by operadic derivations, which is not compatible with the differentials. However, it turns out that the sum of the two actions considered above respects the differentials, and yields the desired action.

%\begin{rem}
Finally, we note that this action does not readily descend to an action of $\fdGC$ on $\fBVGraphs$, the reason being that the right action $\bullet$ does not respect the operadic ideal $I$ considered above. Inserting at a vertex with a tadpole might remove the tadpole.
However, the action descends to an action of the divergence free part $\fdGC_{div}\subset \fdGC$.

\section{Proof of Theorem \ref{thm:main}}\label{sec:theproof}

\subsection{An auxiliary map}
\begin{figure}
\begin{align*}
L&=
\begin{tikzpicture}[baseline={(current bounding box.center)},every edge/.style={draw, -triangle 60}]
%\node at (2,5) {$L=[d,\iota]$};
%\node at (2,.7) {$-$};
%\draw (0,0) ellipse (2cm and 1cm);
%\draw (0,3) ellipse (2cm and 1cm);
%\draw (-2,0)--(-2,3) (2,0)--(2,3);
\node [draw, inner sep=1] (out) at ($(0,0)+(-90:2 and 1)$) {$\vin$};
\node [draw, inner sep=1] (in) at ($(0,3)+(-90:2 and 1)$) {$\vout$};
\node[int] (e1) at (0,.6) {};
\draw[-triangle 60] (e1)--(out) ;
\end{tikzpicture}
-
\begin{tikzpicture}[baseline={(current bounding box.center)},every edge/.style={draw, -triangle 60}]
\node [draw, inner sep=1] (out) at ($(0,0)+(-90:2 and 1)$) {$\vin$};
\node [draw, inner sep=1] (in) at ($(0,3)+(-90:2 and 1)$) {$\vout$};
\node[int] (e1) at (0,0.6) {};
\draw[triangle 60-] (e1)--(in) ;
\end{tikzpicture}
&
F(L)&=
\begin{tikzpicture}[baseline=-.65ex]
 \node[ext] (v) at (0,0) {1};
 \node[int] (w) at (0,1) {};
\draw[-triangle 60] (v) edge (w);
\end{tikzpicture}
+
\begin{tikzpicture}[baseline=-.65ex]
 \node[ext] (v) at (0,0) {1};
\node[int] (w) at (0,1) {};
\draw[-triangle 60] (w) edge (v);
\end{tikzpicture}
+
\begin{tikzpicture}[baseline=-.65ex]
\clip (-.4,-.5) rectangle (.4,2);
 \node[ext] (v) at (0,0) {1};
\node[int] (w) at (0,1) {};
\draw[-triangle 60] (w) edge[loop, looseness=20] (w);
\end{tikzpicture}
\end{align*}
 \caption{\label{fig:Limage}The element $L\in \fdGC_1$ (left) and its image in $\fBVGraphs(1)$ under the map $F$ of Proposition \ref{prop:F} (right).}
\end{figure}
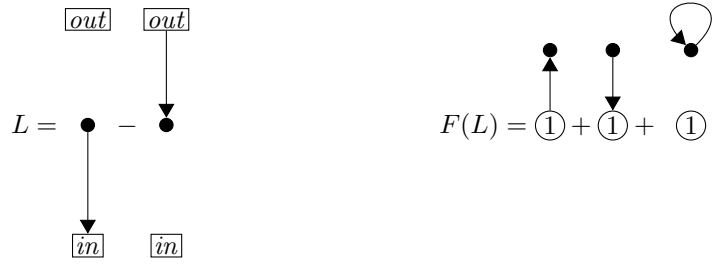

The key step to the proof of Theorem \ref{thm:main} will be the following result, which allows us to identify $\fdGC_1$ with a complex whose cohomology we can compute.
\begin{prop}
\label{prop:F}
There is an isomorphism of complexes
\[
F \colon \fdGC_1 \to \fBVGraphs(1) %\fdGraphs^\whl(1)
\]
such that:
\begin{enumerate}
\item $F$ is compatible with the operadic compositions, i.e., 
\[
F(\Gamma_1\circ \Gamma_1')=F(\Gamma_1)\circ F(\Gamma_1')
\]
for all $\Gamma_1, \Gamma_1'\in \fdGC_1$, where ``$\circ$'' denotes the operadic composition.
\item $F$ is compatible with the right $\fdGC_{div}$ action, i.e., 
\[
F(\Gamma_1\bullet \Gamma)=F(\Gamma_1)\bullet \Gamma
\] 
for all $\Gamma_1\in \fdGC_1$ and $\Gamma\in \fdGC_{div}$.
%\item $F$ is compatible with the differential, i.~e., it is an isomorphism of complexes.
%\item It sends the subspace of $\fdGC_1^\whl$ spanned by graphs containing at least one tadpole to the subspace of $\fdGraphs^\whl(1)$ spanned by graphs containing at least one tadpole at an internal vertex.
\item The image of the element $L\in \fdGC_1$ under $F$ is the element $F(L)$ depicted in Figure \ref{fig:Limage}. 
\end{enumerate}
\end{prop}

% \begin{cor}
% The map $F$ from Lemma \ref{lem:F} descends to a map
% \[
% \bar F \colon \fdGC_1 \to \fBVGraphs(1)'
% \]
% between the quotients $\fdGC_1$ and $\fdGraphs(1)'$ of $\fdGC_1^\whl$ and $\fdGraphs^\whl(1)$ by the graphs with tadpoles at internal vertices. This map $\bar F$ is furthermore a map of complexes, i.e., compatible with the differentials.
% \end{cor}

The proof will occupy the remainder of this subsection.
Let $\Gamma_1$ be a graph in $\fdGC_1$. We define 
\[
F(\Gamma_1) = (-1)^{n_{\vout}} \sum_\Gamma \Gamma \in \fBVGraphs(1)
\]
where $n_{\vout}$ is the valence of $\vout$ in $\fdGC_1$ and the sum runs over all graphs obtained from $\Gamma_1$ by (i) removing vertex $\vout$, (ii) reconnecting the edges previously incident to $\vout$ in some way to vertices in $\Gamma_1$ and (iii) renaming vertex $\vin$ to vertex $1$.
\begin{ex}
To give a concrete example:
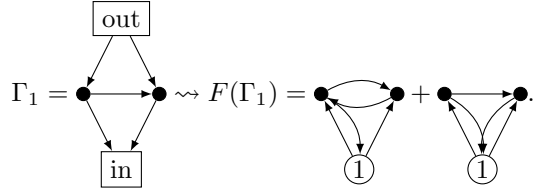
\[
\Gamma_1=
\begin{tikzpicture}[every edge/.style={draw, -latex}, baseline=-.65ex]
 \node[draw] (in) at (0,1) {out}; 
 \node[draw] (out) at (0,-1) {in};
 \node[int] (v1) at (-.5,0) {};
 \node[int] (v2) at (.5,0) {};
\draw (in) edge (v1) edge (v2) (v1) edge (v2) edge (out) (v2) edge (out);
\end{tikzpicture}
\rightsquigarrow
F(\Gamma_1)=
\begin{tikzpicture}[every edge/.style={draw, -latex}, baseline=-.65ex]
 \node[ext] (in) at (0,-1) {1}; 
 \node[int] (v1) at (-.5,0) {};
 \node[int] (v2) at (.5,0) {};
\draw (in) edge (v1) edge (v2) (v1) edge[bend left] (v2) edge [bend left] (in) (v2) edge[bend left] (v1);
\end{tikzpicture}
+
\begin{tikzpicture}[every edge/.style={draw, -latex}, baseline=-.65ex]
 \node[ext] (in) at (0,-1) {1}; 
 \node[int] (v1) at (-.5,0) {};
 \node[int] (v2) at (.5,0) {};
\draw (in) edge (v1) edge (v2) (v1) edge (v2) edge [bend left] (in) (v2) edge[bend right] (in);
\end{tikzpicture}
.
\]
Note that graphs with double edges that occur in the sum are zero and are not shown here. Similarly, by definition of $\fBVGraphs$ graphs with tadpoles at internal vertices are zero.
\end{ex}

We have to show the assertions made in Proposition \ref{prop:F}.
First, we claim that $F$ is an isomorphism of graded vector spaces. To see this filter $\fdGC_1$ by the valence of $\vout$ and filter $\fBVGraphs(1)$ by the number of outgoing edges at the external vertex $1$. Is is easy to see that $F$ is compatible with these filtrations. Consider the associated graded $\gr F$.
It acts on a graph $\Gamma_1\in\fdGC_1$ (up to sign) by (i) connecting vertices $\vin$ and $\vout$ and (ii) renaming the newly formed vertex $1$. 
This is clearly an isomorphism, the inverse map just splits vertex 1 appropriately into $\vin$ and $\vout$.
\begin{rem}\label{rem:explinverse}
From this description one can also easily find an explicit formula for the inverse $F^{-1}$ of $F$. We leave it to the reader. One may use that for an invertible linear map $M=D+N$ with $N$ nilpotent $M^{-1}=D^{-1}- D^{-1}ND^{-1}+D^{-1}ND^{-1}ND^{-1}-\dots$.
\end{rem}

Next, statement (3) of the proposition is an easy explicit computation.
It is also easy to convince oneself that statement (2) is correct. The most difficult assertion is statement (1). 
Fix graphs $\Gamma_1, \Gamma_1'\in\fdGC_1$. We want to show that 
\[
F(\Gamma_1\circ \Gamma_1')=F(\Gamma_1)\circ F(\Gamma_1').
\]
Let us depict $\Gamma_1$ and $\Gamma_1'$ schematically as 
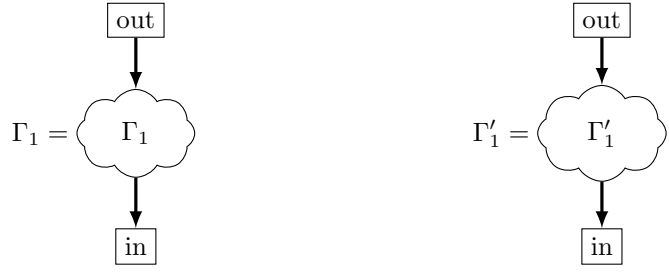
\begin{align*}
\Gamma_1 &=
\begin{tikzpicture}[every edge/.style={draw, very thick, -latex}, baseline=-.65ex]
\node [cloud, draw,cloud puffs=8,cloud puff arc=120, aspect=1.6, inner ysep=.5em] (v1) at (0,0) {$\Gamma_1$};
\node[draw] (out) at (0,1.5) {out}; 
 \node[draw] (in) at (0,-1.5) {in};
\draw (out) edge (v1) (v1) edge (in); 
\end{tikzpicture}
&
\Gamma_1' &=
\begin{tikzpicture}[every edge/.style={draw, very thick, -latex}, baseline=-.65ex]
\node [cloud, draw,cloud puffs=8,cloud puff arc=120, aspect=1.6, inner ysep=.5em] (v1) at (0,0) {$\Gamma_1'$};
\node[draw] (out) at (0,1.5) {out}; 
 \node[draw] (in) at (0,-1.5) {in};
\draw (out) edge (v1) (v1) edge (in); 
\end{tikzpicture}
\end{align*}
where the thick arrows shall stand for (possibly) multiple arrows connecting vertices of $\Gamma_1$, $\Gamma_1'$ to $\vin$ and $\vout$.
A general term (graph) in $\Gamma_1\circ \Gamma_1'$ can be depicted as follows
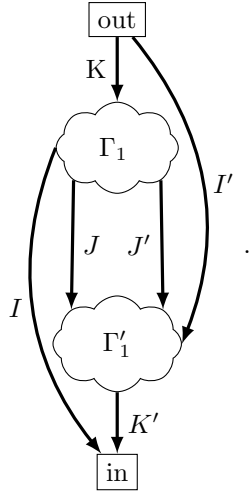
\[
\begin{tikzpicture}[every edge/.style={draw, very thick, -latex}, baseline=-.65ex]
\node [cloud, draw,cloud puffs=8,cloud puff arc=120, aspect=1.6, inner ysep=.5em] (v1) at (0,-1.3) {$\Gamma_1'$};
\node [cloud, draw,cloud puffs=8,cloud puff arc=120, aspect=1.6, inner ysep=.5em] (v2) at (0,1.3) {$\Gamma_1\phantom{'}$}; 
\node[draw] (in) at (0,3) {out}; 
 \node[draw] (out) at (0,-3) {in};
\draw (v1.south) edge node[right] {$K'$} (out) 
 (v2.south east) edge node[left] {$J'$} (v1.north east)
(v2.south west) edge node[right] {$J$} (v1.north west) (v2.west) edge[bend right] node[left] {$I$} (out) (in) edge node[left] {K} (v2) edge[bend left] node[right] {$I'$} (v1.east); 
\end{tikzpicture}
.
\]
Here some subset $I$ of the edges incident at $\vin$ on $\Gamma_1$ is connected to $\vin$, while the remainder $J$ of the edges incident at $\vin$ on $\Gamma_1$ is connected to $\Gamma_1'$ in some way. Similarly, some subset $I'$ of the edges incident at $\vout$ on $\Gamma_1'$ is connected to $\vout$, while the remainder $J'$ of the edges incident at $\vout$ on $\Gamma_1'$ is connected to $\Gamma_1$ in some way. Next consider $F(\Gamma_1\circ \Gamma_1')$. A general term may be depicted as follows:
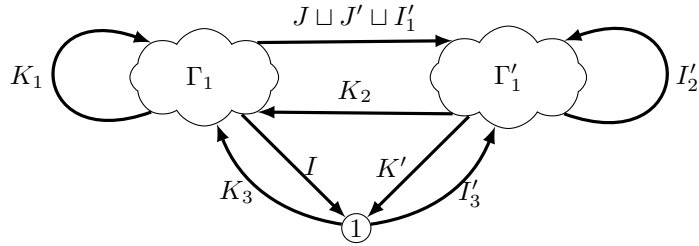
\[
\begin{tikzpicture}[every edge/.style={draw,very thick, -latex}, baseline=-.65ex]
\node [cloud, draw,cloud puffs=8,cloud puff arc=120, aspect=2, inner ysep=.5em] (v1) at (-2,2) {$\Gamma_1\phantom{'}$};
\node [cloud, draw,cloud puffs=8,cloud puff arc=120, aspect=2, inner ysep=.5em] (v2) at (2,2) {$\Gamma_1'$};
\node[ext] (one) at (0,0) {1};
\draw (v1.north east) edge node[above]{$J\sqcup J'\sqcup I_1'$} (v2.north west)
 (v2.south west) edge node[above]{$K_2$} (v1.south east)
(v1) edge node[right] {$I$} (one) 
(v2) edge node[left] {$K'$} (one) 
(v1.south west) edge[out=200, in=160, looseness=5] node[left] {$K_1$} (v1.north west) 
(v2.south east) edge[out=-20, in=20, looseness=5] node[right] {$I_2'$} (v2.north east)
(one) edge[bend left] node[left] {$K_3$} (v1) 
      edge[bend right] node[right] {$I_3'$} (v2); 
\end{tikzpicture}
.
\]
Edges $K$ previously connected to $\vout$ of $\Gamma_1$ are connected either to vertices of $\Gamma_1$, $\Gamma_1'$ or to vertex $1$. 
We split accordingly $K=K_1\sqcup K_2\sqcup K_3$.
Similarly the subset $I'$ of edges as before is further split into $I'=I'_1\sqcup I'_2\sqcup I'_3$, with edges in $I'_1$ being connected to vertices of $\Gamma_1$, while edges in $I'_2$ are connected to vertices of $\Gamma_1'$ and edges in $I'_3$ are connected to vertex $1$. Note that the overall sign of such terms is 
\[
(-1)^{|K|+ |I'|}=(-1)^{|K|+ |I_1'|+|I_2'|+|I_3'|}.
\]
It follows that all terms for which $J'\sqcup I_1'\neq \emptyset$ cancel out. This is because the same edge may participate in either $J'$ or $I_1'$, with opposite signs. Hence we are left with terms for which $I'=J_1'=\emptyset$. 
But these terms are exactly those appearing in $F(\Gamma_1)\circ F(\Gamma_1')$.
Note that in particular that there are all edges in $F(\Gamma_1)\circ F(\Gamma_1')$ between a vertex in $\Gamma_1$ and a vertex in $\Gamma_1'$ originate from edges in $\Gamma_1$ (and not from edges in $\Gamma_1'$).   

We have thus shown that $F$ is an isomorphism of graded vector spaces and that assertions (1)-(3) of Proposition \ref{prop:F} hold. 
It remains to be shown that $F$ is an isomorphism of complexes, i.~e., that it is compatible with the differential.
For $\Gamma_1\in \fdGC_1$ we want to show that 
\[
\delta F(\Gamma_1) = F(\delta \Gamma_1).
\]
Unraveling the formulas for the differentials this can be seen to be equivalent to
\[
(-1)^{|\Gamma_1|} F(\Gamma_1) \bullet \mu
+(-1)^{|\Gamma_1|} F(\Gamma_1) \circ \mu_1 %(\alpha_1+\alpha_2)
+\mu_1 \circ F(\Gamma_1)
=
F(
(-1)^{|\Gamma_1|} \Gamma_1 \bullet \mu
+(-1)^{|\Gamma_1|} \Gamma_1 \circ L
+L \circ \Gamma_1
)
\]
where 
\begin{align*}
 \mu
&=
\begin{tikzpicture}[baseline=-.65ex]
 \node[int](v) at (0,0) {};
\node[int](w) at (0.7,0) {};
\draw[-triangle 60] (v) edge (w);
\end{tikzpicture}
&
\mu_1
&=
\begin{tikzpicture}[baseline=-.65ex]
 \node[ext](v) at (0,0) {1};
\node[int](w) at (0.7,0) {};
\draw[-triangle 60] (v) edge (w);
\end{tikzpicture}
+
\begin{tikzpicture}[baseline=-.65ex]
 \node[int](v) at (0,0) {};
\node[ext](w) at (0.7,0) {1};
\draw[-triangle 60] (v) edge (w);
\end{tikzpicture}.
\end{align*}

By using the second assertion of Proposition \ref{prop:F} (that we already showed) and that $\nabla\mu=0$, we see that the first terms on both sides are the same. By using the first assertion of the proposition it follows that 
\[
F(\Gamma_1 \circ L) = F(\Gamma_1) \circ F(L).
\]
The element $F(L)$ is depicted in Figure \ref{fig:Limage}. Hence
\[
F(\Gamma_1 \circ L) = F(\Gamma_1) \circ F(L)
= F(\Gamma_1) \circ \mu.
\]
Next compute 
\[
F(L\circ \Gamma_1) = F(L)\circ F(\Gamma_1)=\mu\circ F(\Gamma_1).
\]
Hence the above equation holds and we have shown Proposition \ref{prop:F}.
\hfil \qed

% Again the tadpole part vanishes upon passing to quotients, so 
% \[
% \bar F(L\circ \Gamma_1) = \mu\circ \bar F(\Gamma_1).
% \]

%\subsection{Construction of $F$ and proof of Lemma \ref{lem:F}}

\begin{rem}
The motivation behind the definition of $F$ is the following. 
Let $\Tpoly$ be the multivector fields on $\R^n$, %with differential twisted by some Poisson structure $\pi$, 
and let and $\Omega_\bullet$ be the differential forms. These spaces are isomorphic, up to degree shift, with the isomorphism given by sending $\gamma\in \Tpoly$ to the differential form 
\[
\iota_\gamma \omega
\]
where $\omega=dx_1\dots dx_n$ is the standard volume form. The isomorphism $\Tpoly\to \Omega_\bullet[-n]$ induces an isomorphism between the algebras of endomorphisms $\End(\Tpoly)\cong \End(\Omega_\bullet)$.
Now $\fdGC_1$ may be viewed as a graphical version of $\End(\Omega_\bullet)$, while $\fdGraphs(1)$ can be seen as a graphical version of $\End(\Tpoly)$. The map $F$ defined above is just the graphical version of the identification $\End(\Tpoly)\cong \End(\Omega_\bullet)$.
\end{rem}

\subsection{Remainder of the proof of Theorem \ref{thm:main}}
By Proposition \ref{prop:F} we know that 
\[
H(\fdGC_1) \cong H(\fBVGraphs(1)).
\]
The right hand side is known by Proposition \ref{prop:fBVcohom} to be 
\[
H(\fBVGraphs(1))\cong (\K \bbo \oplus \K D)\otimes S(H(\fcGC)[-2])
\]
where $\bbo$ and $D$ correspond to the graphs
\begin{align*}
\bbo &=
\begin{tikzpicture}[baseline=-.65ex]
\node[ext] (v) at (0,0) {$1$};
\end{tikzpicture}
&
D&=
\begin{tikzpicture}[baseline=-.65ex]
\node[ext] (v) at (0,0) {$1$};
\draw (v) edge[loop above] (v);
\end{tikzpicture}.
\end{align*}

Note that $\bbo$ and $D$ are the image of the elements $\bbo\in \fdGC_1$ and $B\in \fdGC_1$ (see \eqref{equ:oneandB}) under the map $F$.

Now we can compute $H(\fEGC)$. Note that as graded vector spaces $\fEGC\cong \fGC\oplus \fGC_1$. The parts of the differential are as follows
\[
 \begin{tikzpicture}
  \node (v) at (0,0) {$\fGC$};
  \node (w) at (3,0) {$\fGC_1$};
  \draw[->] (v) edge[loop, looseness=5] node[above] {$\delta$} (v)
            edge node[above] {$\delta_1=L\bullet \cdot$} (w)
        (w) edge[loop, looseness=5] node[above] {$\delta$} (w);
 \end{tikzpicture}.
\]
We may consider a very simple spectral convergent sequence whose first page sees only the differentials $\delta$, whose second page sees $\delta_1$ and for which all higher differentials vanish.
Taking the cohomology with respect to $\delta$ we obtain 
\[
 H(\fGC) \oplus \K \bbo \oplus \K B \oplus \K \bbo \otimes H(\fGC)[-2] \oplus \K B \otimes H(\fGC)[-2]
\]
by Propositions \ref{prop:F} and \ref{prop:fBVcohom}.
Next, applying $\delta_1$ to a given graph cocyle $\Gamma\in \fGC$ we obtain the linear combination
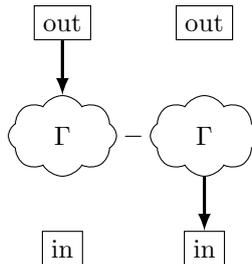
\[
\begin{tikzpicture}[every edge/.style={draw, very thick, -latex}, baseline=-.65ex]
\node [cloud, draw,cloud puffs=8,cloud puff arc=120, aspect=1.6, inner ysep=.5em] (v1) at (0,0) {$\Gamma$};
\node[draw] (out) at (0,1.5) {out}; 
 \node[draw] (in) at (0,-1.5) {in};
\draw (out) edge (v1); 
\end{tikzpicture}
-
\begin{tikzpicture}[every edge/.style={draw, very thick, -latex}, baseline=-.65ex]
\node [cloud, draw,cloud puffs=8,cloud puff arc=120, aspect=1.6, inner ysep=.5em] (v1) at (0,0) {$\Gamma$};
\node[draw] (out) at (0,1.5) {out}; 
 \node[draw] (in) at (0,-1.5) {in};
\draw (v1) edge (in);.
\end{tikzpicture}
\]
To determine its cohomology class in $H(\fGC_1)$ let us apply the map $F$ of Proposition \ref{prop:F}.
We obtain 
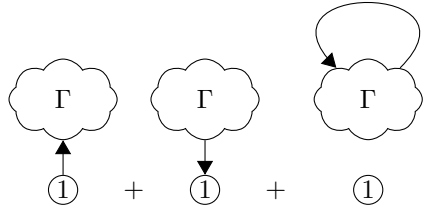
\[
 \begin{tikzpicture}[baseline=-.65ex]
  \node[ext] (v) at (0,0) {1};
\node [cloud, draw,cloud puffs=8,cloud puff arc=120, aspect=1.6, inner ysep=.5em] (v1) at (0,1.2) {$\Gamma$};
\draw[-triangle 60] (v) edge (v1);
 \end{tikzpicture}
+
 \begin{tikzpicture}[baseline=-.65ex]
  \node[ext] (v) at (0,0) {1};
\node [cloud, draw,cloud puffs=8,cloud puff arc=120, aspect=1.6, inner ysep=.5em] (v1) at (0,1.2) {$\Gamma$};
\draw[triangle 60-] (v) edge (v1);
 \end{tikzpicture}
+
 \begin{tikzpicture}[baseline=-.65ex]
\clip (-1,-.5) rectangle (1,3);
  \node[ext] (v) at (0,0) {1};
\node [cloud, draw,cloud puffs=8,cloud puff arc=120, aspect=1.6, inner ysep=.5em] (v1) at (0,1.2) {$\Gamma$};
\draw[-triangle 60] (v1) edge[loop, looseness=5] (v1);
 \end{tikzpicture}
\]
The first two terms together are exact. The last term determines a possibly nontrivial cohomology class, determined by the divergence $\nabla \Gamma$ of $\Gamma$.
Hence the first part of Theorem \ref{thm:main} follows.
To check that indeed $H^0(\fEGC)\cong \grt$, note that $H^{<0}(\fGC)=0$ by \cite{grt}, and hence the divergence operator $\nabla$ has to vanish on $H^0(\fGC)$.
The result then follows since except for $H^0(\fGC)$ no other term contributes to the zero-th cohomology of $\fEGC$.

\section{The explicit formula for divergence free cocycles}\label{sec:divexplicit}
In this section we describe explicitly the cocycles in $\fEGC$ corresponding to divergence free cocyles in $\fGC$.
Concretely, we will describe a map of dg Lie algebras
\[
 \Psi \colon \fGC_{div} \to \fEGC.
\]
In particular, if one picks divergence free representatives of the elements in $\grt\cong H^0(\fGC)$, one obtains an explicit formula for the corresponding cocycles in $\fEGC$.
% \begin{lemma}
% \label{lem:bullet}
% For $\Gamma \in \fGC_1$ and $X\in \fdGC$ divergence free, i.e., 
% \[
% \co{\circlearrowleft
% }{X}=0
% \]
% we have that 
% \[
% F(\Gamma \bullet X) = F(\Gamma) \bullet X.
% \]
% \end{lemma}
% \begin{proof}
% 
% \end{proof}

% We need the following preparatory Lemma.
% \begin{lemma}
% \label{lem:circ}
% Let $\Gamma, \Gamma' \in \fGC_1$ be graphs such that there is no (direct) edge from vertex $\vout$ to vertex $\vin$. Then
% \[
% F(\Gamma \circ \Gamma') =  F(\Gamma) \circ F(\Gamma').
% \]
% \end{lemma}
% \begin{proof}
% 
% \end{proof}

For a element $X \in \fGC_{div}$, let $X_1\in \fdGraphs(1)$ be the element obtained by declaring one vertex (say the first) external.
Then we define the map 
\begin{gather*}
\Psi\colon \fGC_{div} \to \fEGC \\
X \mapsto X + F^{-1}(X_1).
\end{gather*}
We claim that the map is a map of dg Lie algebras.
To check compatibility with the differential, we have to verify that 
\[
F^{-1}((\delta X)_1) =  L\bullet X + \delta F^{-1}(X_1).
\]

Since $F$ is an isomorphism we may as well apply $F$ to the both sides and check that
\[
(\delta X)_1 = F(L\bullet X + \delta  F^{-1}(X_1))
= F(L\bullet X)+ \delta X_1
\]
Now consider the first term on the right. Since $X$ is divergence free we have by Proposition \ref{prop:F} that
\[
 F(L\bullet X) =  F(L) \bullet X = \mu_1 \bullet X.
\]
%where $\pi \fdGraphs(1)^\whl \to \fdGraphs(1)'$ is the projection to the quotient. We computed $F(L)$ above. Since $X$ is divergence free it follows that 
%\[
%\pi ( F(L)\bullet X) = \mu_1 \bullet X.
%\]
But it is not hard to check by a small graphical calculation that 
\[
 \mu_1 \bullet X + \delta X_1 = (\delta X)_1
\]
and hence compatibility with the differential follows.

Next let us consider compatibilty with the Lie bracket. Let $Y$ be another closed element in $\fGC_{div}$ and let again $Y_1$ be the element in $\fdGraphs(1)$ obtained by declaring one vertex external. We assume $X$ and $Y$ are homogeneous of degrees $|X$ and $|Y$.
%Note that $F^{-1}(X_1)$ and $F^{-1}(Y_1)$ do not contain graphs with edges $\vout \to \vin$.
Compute
\begin{multline*}
\co{X+F^{-1}(X_1) }{Y+F^{-1}(Y_1)}
=\\ \co{X}{Y}
+
F^{-1}(X_1)\bullet Y
-(-1)^{|X||Y|}
F^{-1}(Y_1) \bullet X
+
 F^{-1}(X_1)\circ F^{-1}(Y_1)
-
(-1)^{|X||Y|}
F^{-1}(Y_1)\circ F^{-1}(X_1).
\end{multline*}
Apply $F$ to the part in $\fdGC_1$ and use Propostion \ref{prop:F} again to compute

\begin{align*}
F(\dots)
&= 
 F( F^{-1}(X_1)\bullet Y))
-
(-1)^{|X||Y|}
 F( F^{-1}(Y_1) \bullet X)
+
 F( F^{-1}(X_1)\circ  F^{-1}(Y_1))
\\&\quad\quad\quad\quad\quad\quad\quad\quad\quad\quad\quad\quad -
(-1)^{|X||Y|}
 F( F^{-1}(Y_1)\circ  F^{-1}(X_1))
\\
&= X_1\bullet Y
-
(-1)^{|X||Y|}
Y_1 \bullet X
+
X_1\circ Y_1
-
(-1)^{|X||Y|}
Y_1\circ X_1
.
\end{align*}
But this is $(\co{X}{Y})_1$ and we are done. %(TODO: show or add ref)

%Since graphs with tadpoles do not span an ideal under compositions, the tadpole part of $F(L)$ may now contribute. 

\begin{rem}
 Note that the images of $\Psi$ commute with the element $B$ of eqn. \eqref{equ:oneandB}.
This means in particular that the derivation of $\Omega_\bullet$ corresponding to elements of $\grt$ is compatible with the de Rham differential.
Hence the main result of \cite{mech} of the compatibility of the Shoikhet morphism with the de Rham differential may be extended to all formality morphisms obtained by the $\fGC_{div}$ action.
\end{rem}

%\nocite{*}
\bibliographystyle{plain}
\bibliography{../biblio}

\end{document}